\ifx\LocalElevenPtArticleFormatLoaded\undefined
  \documentclass[a4paper,11pt]{article}     
  \usepackage{array}                        
  \usepackage{theorem}                      
  \usepackage{amsmath,amscd,amssymb}        
  \usepackage{latexsym}                     
  \usepackage[german,english]{babel}        
  \usepackage[applemac]{inputenc}
 \usepackage{xypic}
\fi

\newtheorem{theorem}{Theorem}[section]
\usepackage{theorem}

\newtheorem{corollary}[theorem]{Corollary}
\newtheorem{proposition}[theorem]{Proposition}
\theorembodyfont{\rm}

\newcommand{\CC}{{\mathbb{C}}}

\newcommand{\PP}{{\mathbb{P}}}

\newcommand{\ZZ}{{\mathbb{Z}}}

\newcommand{\DD}{{\mathbb{D}}}

\newenvironment{Proof}{\begin{ProofwCaption}{Proof}}{\end{ProofwCaption}}
\newenvironment{Proof*}[1]{\begin{ProofwCaption}{{#1}}}{\end{ProofwCaption}}
\newenvironment{ProofwCaption}[1]%
  {\addvspace\theorempreskipamount \noindent{\it #1.}\rm}%
  {\qed \par \addvspace\theorempostskipamount}
\newcommand{\qedsymbol}{\mbox{$\Box$}}
\newcommand{\qed}{\quad\qedsymbol}
\begin{document}
\begin{center}
{\bf\Large A bound on the irregularity of abelian scrolls in projective
space}\\[3mm]
{\large C. Ciliberto and K. Hulek}
\end{center}
\begin{center}
{\em{Dedicated to Professor Hans Grauert}}
\end{center}

\noindent {\bfseries Abstract.}  {\footnotesize
We prove that the irregularity of a smooth {\em abelian scroll} whose
dimension is
at least half of that of the surrounding projective space is bounded by $2$.
We also discuss some existence results and open problems.
}

\section{Introduction}
Let $A\subset\PP^{l-1}$ be an $n$-dimensional abelian variety
and let $ G\subset A$ be a subgroup of
order
$k$. Then one can define for every point $P\in A$ the linear space
$$
S^G(P):=\mathrm {Span}\langle P+\rho ; \rho \in G\rangle.
$$
In general one expects $S^G(P)$ to have dimension $k-1$ and the union
of these spaces
$$
Y=\bigcup_{P\in A} S^G(P)
$$
is then a (possibly singular) scroll of expected dimension $n+k-1$.
The classical example of this construction is the case where $E$ is an
elliptic quintic
normal
curve in $\PP^4$ and $G=\ZZ_2$. This leads to the quintic elliptic scroll,
the only
smooth irregular scroll in $\PP^4$.

We shall refer to a smooth scroll $Y$
constructed
as above and having the expected dimension
as  the {\em abelian scroll} determined by the pair $(A,G)$.
In particular, if $\rho\in A_{(k)}$ is a point of
order
$k$, we can consider the cyclic subgroup $G\simeq \ZZ_k$ generated by $\rho$.
We call the corresponding
abelian scroll the {\em cyclic scroll} determined by the pair $(A,\rho)$.

Observe that for smooth abelian scrolls $Y$ we have
$2\dim Y\le l-1$. This follows immediately from Barth's theorem which says that
any smooth subvariety $Y \subset \PP^{l-1}$ with $2 \dim Y > l-1$
has irregularity $q(Y)=0$. The
main result of this paper is the following:

\begin{theorem}\label{theo0.1}
Let $Y\subset \PP^{l-1}$ be a smooth abelian scroll with $2\dim Y=l-1$. Then
$q(Y)\le 2$, i.e. $A$ is either an elliptic curve or an abelian surface.
Moreover, in
this case $Y$ (or equivalently $A$) is linearly normally embedded.
\end{theorem}
This non-existence result can be seen as an extension of a result
of Van de Ven \cite{V} which says that an abelian variety $A$ of dimension
$n$ can be embedded in $\PP^{2n}$ only if $n=1$ or $2$. We shall briefly
discuss the existence of abelian scrolls and some open problems.

As we
indicate at the end of section 2,  the question considered here is related to
the general interesting
problem of classifying  smooth
projective varieties whose dimension equals their codimension, a problem
motivated by questions of
Griffiths and van de Ven concerning smooth surfaces in $\PP^4$.

{\footnotesize
\noindent {\bf Acknowledgement.}
We warmly thank the
referee for carefully reading the paper and, in particular, for having
suggested an elegant
simplification of the first version of this note which enabled us to avoid
the use of some binomial
identities. However, we thank M.~Ern\'{e} and S.~Rao for having given us the
right references for
these identities in the first place.
Both authors gratefully acknowledge partial support by the EU project
EAGER
(HPRN-CT-2000-00099). The second author was also partially supported by the
DFG Schwerpunktprogramm
``Globale Methoden in der komplexen Geometrie''. }

\section{Non-existence of scrolls}

This section is devoted to the proof of the non-existence
theorem (\ref{theo0.1}).
The main
ingredients in the proof are the double point formula and some estimates
involving binomial coefficients.

Before we can prove our result we need some preparations. Consider the \'etale
$k:1$ quotient
$$
\bar{\pi}:A\rightarrow\bar{A}:=A/ G
$$
and let ${\cal L}$ be the line bundle on $A$ which defines the embedding into
$\PP^{l-1}$. Then
$$
{\cal E}:=\bar{\pi}_{\ast}{\cal L}
$$
is a rank $k$ vector bundle on $\bar{A}$ and we set
$$
X:=\PP({\cal E}).
$$
The natural map $\bar{\pi}^{\ast}{\cal E} \to {\cal L}$ induces an
inclusion of $A$ into $X$ as a multi-section and
the linear system $|{\cal O}_X(1)|=|{\cal O}_{\PP({\cal E})}(1)|$
induces the complete linear system $|{\cal L}|$ on $A$ and
maps $X$
onto a scroll $Z$ (see \cite [Lemma 2.1] {CH}).
If $Y$ is embedded by a complete
linear system, then $Z=Y$, otherwise $Y$ is a
projection of $Z$. Also as in \cite[p. 360]{CH} we have
$$
\bar{\pi}^{\ast}{\cal E}={\cal L}\oplus t^{\ast}_{\rho_1}{\cal
L}\oplus\ldots\oplus t^{\ast}_{\rho_{k-1}}{\cal L}
$$
where $G=\{1,\rho_1,\ldots,\rho_{k-1}\}$.
Note that topologically $\tilde{X}:=\PP(\bar{\pi}^{\ast}{\cal E})
\cong{A}\times\PP^{k-1}$ is trivial. We can also choose another covering
$\hat{\pi}: \hat{A} \to A$ of order $h$ where $h$ divides $k$ such that
$\hat{X}:=\PP(\hat{\pi}^{\ast}\bar{\pi}^{\ast}{\cal
E})\cong\hat{A}\times\PP^{k-1}$ as algebraic varieties.

{\it Proof of theorem 0.1}
Let $m=\dim Y=n+k-1$. Then the double point formula
(see \cite [Theorem (9.9.3)]{Fu}) reads
$$
[\DD]=\varphi^{\ast}\varphi_{\ast}[X]-c_m(N_{\varphi})\cap[X]
$$
where $\varphi:X\rightarrow Y$ is the map given by $|{\cal O}_X(1)|$,
possibly followed by a projection. We first
determine the contribution $\varphi^{\ast}\varphi_{\ast}[X]$. We have
already observed that topologically
$\tilde{X}$ is a product and we can write for the hyperplane section on
$\tilde{X}$:
$$
\tilde{H}=c+h
$$
where $c=c_1({\cal L})\in H^2(A,\ZZ)$ and $h$ is the positive generator of
$H^2(\PP^{k-1},\ZZ)$. Hence
$$
\tilde{H}^m=\tilde{H}^{n+k-1}=(c+h)^{n+k-1}={\binom {n+k-1} {k-1}} c^n.
$$
Since $\tilde{X}\rightarrow X$ has degree $k$ we find that
$$
\varphi^{\ast}\varphi_{\ast}[X]=\frac 1{k^2}{\binom{n+k-1} {k-1}}^2 (c^n)^2.
$$
Instead of the normal bundle of $\varphi$ we compute the normal bundle of
$\tilde{\varphi}:=\varphi\circ f$ where $f:\tilde{X}\rightarrow X$ is the
\'etale $k:1$ map induced from $\bar{\pi}:A\rightarrow \bar{A}$.
{}From the exact
sequence
$$0\rightarrow T_{\tilde{X}}\rightarrow\tilde{\varphi}^{\ast}
T_{\PP^{l-1}}\rightarrow N_{\tilde{\varphi}}\rightarrow 0
$$
we find
$$
\begin{array}{rcl}
c(N_{\tilde{\varphi}}) &=&
c(\tilde{\varphi}^{\ast} T_{\PP^{l-1}}) c (T_{\tilde{X}})^{-1}\\[1mm]
 &=& (1+c+h)^l(1+h)^{-k}
\end{array}
$$
where we have again used that topologically
$\tilde{X}={A}\times\PP^{k-1}$. Now we notice that $c^{n+1}=h^k=0$ and,
therefore, the top
term of $c(N_{\tilde{\varphi}})$ is just the coefficient of the
monomial $c^nh^{k-1}$, which is the
same as the coefficient of the same monomial in
$${{l}\choose n}c^n(1+h)^{l-n}(1+h)^{-k}={{l}\choose n}c^n(1+h)^{l-n-k}.$$
Since $l=2n+2k-1$ we get
$$c_m(N_{\tilde{\varphi}})=c^n{{n+k-1}\choose n}{{2n+2k-1}\choose n}.$$
Altogether we find for the class of the double locus $\DD$ that
$$
[\DD]=\frac 1 k c^n\left[
\frac 1 k \binom{ n+k-1}
{k-1}^2 c^n-{{n+k-1}\choose n}{{2n+2k-1}\choose n}
\right].
$$
Using Riemann-Roch on the abelian variety $A$ we find
$$
c^n \geq n! (2n+2k-1)
$$
with equality if and only if the linear system is complete.
Then the assertion of the theorem follows if one shows that
$$
 {\binom{ n+k-1} {k-1}} (2n+2k-1) n!\ge
k \binom{ 2n+2k-1} {n}
$$
where equality holds if and only if $n=1$ or $2$.
Indeed, it is a straightforward calculation to check equality
for $n=1$ and $n=2$. Now
assume $n\ge 3$. Evaluating the binomial coefficients and cancelling terms,
the above inequality becomes
$$
n!\prod\limits^n_{l=2} (n+k-l+1)\ge \prod\limits^n_{l=2}(2n+2k-l)
$$
which can be rewritten as
$$
\prod\limits^n_{l=2}(n+k-l+1) l\ge\prod\limits^n_{l=2}(2n+2k-l).
$$
To check this inequality it is enough to check it termwise. The inequality
$$
(n+k-l+1) l\ge 2n+2k-l\quad\mbox{ for } l=2,\ldots, n
$$
is equivalent to
$$
n+k\ge l\qquad \mbox{ for } l=2,\ldots, n
$$
which is trivially true.
\hfill $\Box$

\section{Existence of scrolls}
In this section we shall discuss conditions under which a linearly
normal abelian scroll $Y$ determined by a pair $(A,G)$
is smooth. We shall keep the notation as in section 1 and we shall start
with the scroll $X=\PP({\cal E})$. Recall
that
$A$ is contained in $X$ as a $k$-section and that
${\cal O}_X(1)|_A\cong{\cal L}$,
essentially by the definition of ${\cal O}_X(1)$. We shall always assume that
${\cal L}$ is very ample on $A$. The line bundle ${\cal L}$ is called
$(m-1)$-very ample, if for every cluster
(i.e. $0$-dimensional subscheme) $\zeta\subset A$ of length $\le m$ the
restriction map
$H^{0}(A,{\cal L})\rightarrow H^{0}(\zeta,{\cal L}\otimes{\cal
O}_{\zeta})$ is surjective and hence $\varphi_{|\cal L|}(\zeta)$
spans a $\PP^{m-1}$.

\begin{proposition}\label{prop2.1}
{\rm(i)}\ If ${\cal L}$ is $(k-1)$-very ample, then $|{\cal O}_X(1)|$ is
base point free and $\varphi:X\rightarrow Y$ is a finite map.

\noindent{\rm(ii)} If ${\cal L}$ is $k$-ample, then $\varphi:X\rightarrow Y$ is
birational and an embedding near $A$.
\end{proposition}

\begin{Proof}
(i)\  $A$ intersects each fibre of $X$ in $k$ independent points. This follows
since by construction of $X$ these $k$ points are the image of points in
$\tilde{X}$ given by the sections which come from the splitting of
$\bar{\pi}^{\ast}{\cal E}$. Since these points are mapped to independent
points in
$\PP^{l-1}$ it follows that $|{\cal O}_X(1)|$ restricted to a fibre of $X$
cannot have base points. In order to show that $\varphi$ is finite, it is
enough
to prove that ${\cal O}_X(1)$ has positive degree on each curve $C$ in $X$.
This
can be checked by pulling back to the product
$\hat{X}=\hat{A}\times\PP^{k-1}$ where we have a product polarization.

\noindent (ii)\  It follows immediately from our assumption that every
$\PP^{k-1}$ of the scroll $Y$ meets $A$ only in the $k$ points which span it.
Hence for every point $P \in A$ we have that $\varphi^{-1}(P)$ consists of only
one point. Similarly we prove that $d{\varphi}$ is injective along $A$. Now
assume that $\varphi:X\rightarrow Y$ has degree $d\ge 2$. Over each point of
$A\subset Y$ we have only one preimage. Hence $\varphi$ is ramified along $A$,
contradicting what we have just proved.
\hfill\end{Proof}

\begin{proposition}\label{prop2.2}
If ${\cal L}$ is $(2k-1)$-very ample, then $\varphi:X\rightarrow Y$ is an
isomorphism and, in particular, $Y$ is smooth.
\end{proposition}

\begin{Proof}
The assumption gives immediately that, for two points $P$ and $Q$ whose
difference on $A$ lies in the group
$G$, we have
$S^{k-1}(P)\cap S^{k-1}(Q)=\emptyset$. Hence the map $\varphi$ is injective.
It remains to prove that the
differential
$d{\varphi}$ is injective everywhere. The map $\varphi:X\rightarrow Y \subset
\PP^{l-1}$ embeds every fibre $\pi_x$ of $\pi:X\rightarrow \bar{A}$
through a point $x \in X$ as a
$(k-1)$-dimensional subspace which, by abuse of notation, we also denote by
$\pi_x$. We consider the map
$$
\begin{array}{rcl}
\gamma : A & \rightarrow & G(k-1, l-1)=:\mathrm{Gr}\\
x          &\mapsto      & \pi_x.
\end{array}
$$
where $\pi_x$ is the unique fibre of $X$ containg the point $x$. This map
factors through $\bar{A}$.
Taking the projectivised differential of this map gives us a map
$$
\begin{array}{rcl}
\omega=\PP(d\gamma):\PP(T_{A,x}) & \rightarrow &
\PP(T_{\mathrm{Gr},\pi_x})={\operatorname Pr} (\pi_x,\pi^{0})\\
t &\mapsto & \omega_t
\end{array}
$$
where $\pi^{0}$ is a complementary subspace of $\pi_x$ and
$\operatorname{Pr}(\pi_x,\pi^{0})$ is the projective space of projective
transformations of
$\pi_x$ to $\pi^{0}$. (If $\pi_x=\PP(U)$  and
$\pi^{0}=\PP(W)$ then $\operatorname{Pr}(\pi_x,
\pi^{0})=\PP(\operatorname{Hom}(U,W$)).) Here the identification
$\PP(T_{\mathrm{Gr},\pi_x})=\operatorname{Pr}(\pi_x,\pi^0)$
comes from the fact that
there
is a canonical isomorphism $T_\mathrm{Gr}=\operatorname{Hom}({\cal U},{\cal
V})$ where
${\cal U}$ and ${\cal V}$
are the canonical sub bundle resp. the
quotient bundle. We also have the Gauss map
$$
\begin{array}{rcl}
\Gamma : X &\dashrightarrow& G(n+k-1,l-1)\\
P &\longmapsto & \bar{T}_{X,P}
\end{array}
$$
where $\bar{T}_{X,P}$ is the projective closure of the image of the
differential of $\varphi$ at $P$, considered as a subspace through the point
$P$. The relation between $\gamma$ and $\Gamma$ is the following. Let $P\in X$
and let $x \in \bar{\pi}^{-1}(\pi(P))$. Then
$$
\bar{T}_{X,P}=<\pi_x, \omega_t(P);\ t\in \PP(T_{A,x})>.
$$
In order to check that $d{\varphi}$ is injective at $P$ we have to prove that
the projective map
$$
\begin{array}{rcl}
\omega(P):\ \PP(T_{A,x}) & \dashrightarrow & \pi^0\\
t & \longmapsto & \omega_t (P)
\end{array}
$$
is injective, i.e. is well defined. If this is not the case, then we
have a tangent direction
$t\in\PP(T_{A,x})$ such that the linear map associated to the map
$$
\begin{array}{rcl}
\omega_t:\ \pi_x & \dashrightarrow & \pi^0\\
P & \longmapsto & \omega_t (P)
\end{array}
$$
has a kernel. Assume that this is the case. Consider the germ of a holomorphic
curve $(\CC,0)\rightarrow (A,x)$ which represents the tangent direction $t$ at
$x$. This determines a family of $(\PP^{k-1})'s$ given by
$$
\pi_s=<z_0(s),\ldots,z_{k-1}(s)>.
$$
Choosing suitable coordinates in $\PP^{l-1}$ we can assume that
$$
z_i(0)=e_i\ ;\ i=0,\ldots, k-1
$$
with $e_0,\ldots,e_{l-1}$ the standard basis. For the complementary space
$\pi^0$ we can choose $\pi^0=<e_{k},\ldots,e_{l-1}>$. It is an easy local
computation to check that with respect to the coordinates on $\pi_x$ and
$\pi^0$ given by $e_0,\ldots,e_{k-1}$ and $e_k,\ldots,e_{l-1}$ the linear map
associated to the map
$$
\omega_0:\pi_x\rightarrow\pi^0
$$
is given by the matrix
$$
M=\left(
\begin{array}{cccc}
z'_{0k} & z'_{1k}& \ldots & z'_{k-1,k}\\
\vdots & \vdots & & \vdots\\
z'_{0,l-1} & z'_{1,l-1}& \ldots & z'_{k-1,l-1}
\end{array}
\right) \ (s=0).
$$
Hence $\omega_0$ has a kernel if and only if $\operatorname{rank} M < k$. But
this is equivalent to the assertion that $z_0(0),\ldots, z_{k-1}(0),
z'_0(0),\ldots, z'_{k-1}(0)$ are dependent. In particular there is a cluster of
length $2k$ on
$A$ which contradicts $(2k-1)$-very ampleness of ${\cal L}$.
\hfill\end{Proof}

We can apply this result immediately to elliptic normal curves $E\subset
\PP^{m-1}$. Recall that any $m-1$ points on $E$ (possibly infinitely near) are
linearly independent (otherwise we could find a hyperplane through these
$m-1$ points and {\em{any}} other point on $E$.)

\begin{corollary}\label{cor2.3}
Let $m=2n+1$ be an odd integer. Then there exists an $n$-dimensional elliptic
curve scroll $Y$ in $\PP^{2n}$.
\end{corollary}

\begin{Proof}
Let $E\subset\PP^{m-1}=\PP^{2n}$ be an elliptic normal curve of degree $m$ and
choose a subgroup $G$ of order $n$ of $E$, e.g. a cyclic subgroup.
We want to prove that the
abelian scroll $Y$ determined by the pair $(E,G)$ is smooth.
This follows from Proposition \ref{prop2.2}
since the embedding line bundle ${\cal L}$ is $(m-2=2n-1)$-very ample.
\hfill\end{Proof}

Conversely we can also use our results to bound the very-ampleness of line
bundles on abelian varieties.

\begin{corollary}\label{cor2.3a}
Let ${\cal L}$ be a line bundle on an abelian variety $A$ of dimension $\dim
A=n\ge 3$ with $h^0(A,{\cal L})=l$. Assume that ${\cal L}$ is $k$-very ample
for some odd integer $k$. Then $k<l-2n$.
\end{corollary}

\begin{Proof}
For $k=1$ this is Van de Ven's result \cite{V}. If  $k\ge 3$ we can write
$k=2s-1$ for some $s\ge 2$. The assertion then follows by combining
Theorem (\ref{theo0.1}) and Proposition (\ref{prop2.2}),
applied to the abelian scroll determined by a
pair $(A,G)$, where $G$ can be chosen as any subgroup of order $k$ of $A$,
e.g.
a cyclic subgroup.
{}\hfill\end{Proof}

The crucial open question which remains is the existence of abelian surface
scrolls $Y \subset \PP^{l-1}$ with $2\dim Y=l-1$. In \cite{CH} we constructed
one non-trivial such example, namely the following.\\

\noindent{\bf Example.} Let $A\subset\PP^6$ be a general (i.e.
$\operatorname{rank} \operatorname{NS}(A) =1$) abelian surface embedded
by a line bundle ${\cal L}$ of type $(1,7)$ and let ${\varepsilon}\in A_{(2)}$
be a non-zero 2-torsion point.
Then the cyclic scroll determined by the pair $(A,\varepsilon)$
is smooth of dimension 3 and has degree 21.\\

 The next step would be to
investigate cyclic $\PP^2$-scrolls determined by pairs $(A,\rho)$
where  $A\subset\PP^8$ is an abelian
surfaceof degree 18, embedded by a complete linear system of type $(1,9)$
and $\rho\in A_{(3)}$ is a 3-torsion point. More generally we pose the\\

\noindent{\bf Problem}. Let $A\subset\PP^{2n}$, $n\geq 4$ be a general abelian
surface of degree $4n$ and let $\rho \in A_{(n-1)}$ be a non-zero
$(n-1)$-torsion point. Is the cyclic scroll defined by the pair $(A,\rho)$
smooth?\\

A positive answer to this question would provide us with an infinite series
of smooth scrolls $Y\subset\PP^{2n}$ whose dimension is half that of the
surrounding space and with irregularity 2. Van de Ven has asked the
question whether the
irregularity of smooth surfaces in $\PP^4$ is bounded by 2. As far as we
know, this
question is
still open. In fact, we do not know of any smooth subvariety $Y\subset\PP^{2n}$
with $\dim Y=n\ge 2$ and $q\geq 3$. Let us, therefore, pose the\\

\noindent{\bf Problem}. Give examples of smooth subvarieties
$Y\subset\PP^{2n}$
with $\dim Y=n\geq 2$ and $q\geq 3$.


%
\bibliographystyle{alpha}


\noindent Ciro Ciliberto\par
\noindent Universit\`a di Roma Tor Vergata\par
\noindent Dipartimento di Matematica\par
\noindent Via della Ricerca Scientifica\par
\noindent I 00173 Roma (Italy)\par
\noindent e-mail: Cilibert@axp.mat.uniroma2.it\par

\vskip 0.5truecm

\noindent Klaus Hulek\par
\noindent Institut f\"ur Mathematik\par
\noindent Universit\"at Hannover\par
\noindent Post fach 30060\par
\noindent D 3167 Hannnover (Germany)\par
\noindent e-mail: Hulek@math.uni-hannover.de\par

\end{document}